\documentclass[final,3p,times,12pt]{elsarticle}

\usepackage{lineno}
\usepackage{amsmath}
\usepackage{amssymb}
\usepackage{amsthm}
\usepackage{cases}
\usepackage{url}
\usepackage[pdfpagelabels,plainpages=false]{hyperref}
\hypersetup{ linkcolor =blue,
 citecolor = cyan,
 urlcolor = blue,
colorlinks = true, bookmarksopen = true, bookmarksnumbered = true,
bookmarksopenlevel = {1}, linktocpage = true,
pdfstartview = FitH,
pdftitle = {Mixed Recurrence equations, interlacing zeros, bounds for zeros},
pdfsubject
= {article},
pdfauthor = {Tcheutia et al.},
pdfcreator = {pdftextify}, pdfproducer = {pdftextify}
}

\newtheorem{theorem}{Theorem}
\newtheorem{corollary}[theorem]{Corollary}
\newtheorem{lemma}[theorem]{Lemma}

\newtheorem{remark}[theorem]{Remark}
\def\a{\alpha}
\def\b{\beta}
\def\l{\lambda}
\newcommand{\nn}{\mathbb{N}}

\newcommand{\hypergeom}[5]{\mbox{$
_#1 F_#2\left( \! \left.
\begin{array}{c}
\multicolumn{1}{c}{\begin{array}{c} #3
\end{array}}\\[1mm]
\multicolumn{1}{c}{\begin{array}{c} #4
            \end{array}}\end{array}
\! \right| \displaystyle{#5}\right) $} }
\newcommand{\qhypergeom}[5]{\mbox{$
_#1 \phi_#2\left(\!  \left.
\begin{array}{c}
\multicolumn{1}{c}{\begin{array}{c} #3
\end{array}}\\[1mm]
\multicolumn{1}{c}{\begin{array}{c} #4
            \end{array}}\end{array}
  \right| \displaystyle{#5}\right) $} }

\begin{document}
	
\begin{frontmatter}
\title{Recurrence equations involving different orthogonal polynomial sequences and applications}

\author[label1]{A.S.~Jooste}\ead{alta.jooste@up.ac.za}
\author[label2]{D.D.~Tcheutia}\ead{duvtcheutia@yahoo.fr}
\author[label2]{W.~Koepf\corref{cor1}}\ead{koepf@mathematik.uni-kassel.de}

\ead[url]{www.mathematik.uni-kassel.de/~koepf}
\cortext[cor1]{Corresponding author}
\def\l{\lambda}
\address[label1]{Department of Mathematics and Applied Mathematics, University of Pretoria, Pretoria 0002, South Africa}
\address[label2]{Institute of Mathematics, University of Kassel, Heinrich-Plett Str. 40, 34132 Kassel, Germany}

\begin{abstract}
Consider $\{p_n\}_{n=0}^{\infty}$, a sequence of polynomials orthogonal with respect to $w(x)>0$ on $(a,b)$, and polynomials $\{g_{n,k}\}_{n=0}^{\infty},k \in \nn_0$,
orthogonal with respect to $c_k(x)w(x)>0$ on $(a,b)$, where $c_{k}(x)$ is a polynomial of degree $k$ in $x$. We show how Christoffel's formula can be used to obtain mixed three-term recurrence equations involving the polynomials $p_n$, $p_{n-1}$ and $g_{n-m,k},m\in\{2,3,\dots, n-1\}$. In order for the zeros of $p_n$ and $G_{m-1}g_{n-m,k}$ to interlace (assuming $p_n$ and $g_{n-m,k}$ are co-prime), the coefficient of $p_{n-1}$, namely $G_{m-1}$,  should be of exact degree $m-1$, in which case restrictions on the parameter $k$ are necessary. The zeros of $G_{m-1}$ can be considered to be inner bounds for the extreme zeros of the (classical or $q$-classical) orthogonal polynomial $p_n$ and we give examples to illustrate the accuracy of these bounds. Because of the complexity the mixed three-term recurrence equations in each case, algorithmic tools, mainly Zeilberger's algorithm and its $q$-analogue, are used to obtain them.
\end{abstract}

\begin{keyword}
Classical orthogonal polynomials \sep $q$-classical orthogonal polynomials \sep  mixed three-term recurrence equations \sep interlacing of zeros \sep bounds for zeros
\MSC[2010] 33C05 \sep 33C45 \sep 33F10 \sep 33D15
\end{keyword}

\end{frontmatter}

\section{Introduction }
A sequence $\{p_n\}_{n=0}^{\infty}$ of real polynomials, where $p_n$ is of exact degree $n$, is orthogonal with respect to an absolutely continuous measure that can be represented by a real, positive weight function $w(x)$  on the  (finite or infinite) interval $(a,b),$ if
\begin{equation*}
\int_{a}^{b}p_n(x)p_m(x)w(x)dx= 0,m\neq n.
\end{equation*}
The classical orthogonal  polynomials  considered in this paper  are defined in terms of the generalized hypergeometric series
\[\hypergeom{p}{q}{a_1,\ldots,a_p}{b_1,\ldots, b_q}{x}=\sum_{m=0}^{\infty}\frac{(a_1)_m\cdots (a_p)_m}{(b_1)_m\cdots (b_q)_m}\frac{x^m}{m!},\]
where  $(a)_m$ denotes the Pochhammer symbol (or shifted
factorial)  defined by  \[(a)_m=\left\{
                  \begin{array}{ll}
                    1 & \text{ if } m=0 \\
                    a(a+1)(a+2)\cdots (a+m-1)
                     & \text{ if } m\in \mathbb{N}.
                  \end{array}
                \right.
\]
Their $q$-orthogonal analogues, $0<q<1$, are given in terms of  basic hypergeometric series \cite[Section 1.10]{KLS}
 \[
 \qhypergeom{r}{s}{a_1,\ldots,a_r}{b_1,\ldots,b_s}{q;z}=
 \sum_{k=0}^\infty\frac{(a_1,\ldots,a_r;q)_k}{(b_1,\ldots,b_s;q)_k}\left((-1)^k
 q^{\binom{k}{2}}\right)^{1+s-r}\frac{z^k}{(q;q)_k},
 \]
 where the  $q$-Pochhammer symbol  $(a_1,a_2,\ldots,a_k;q)_n$ is defined by \[(a_1,\ldots,a_r;q)_k:=(a_1;q)_k\cdots(a_r;q)_k,~\text{  with }~
 (a_i;q)_k=\left\{
 \begin{array}{ll}
 \prod\limits_{j=0}^{k-1}(1-a_iq^j)& \text{ if }\ k\in\{1,2,3,\ldots\} \\
  1 &  \text{ if }\  k=0.
  \end{array}
 \right.
 \]
Let the zeros of $p_n$ be $x_{n,1}<x_{n,2}<\dots<x_{n,n}.$ It is well known that $p_n$ satisfies a  three-term recurrence equation
\begin{equation}p_{n}(x)=\label{TTRR1}(x-C_n)p_{n-1}(x)-\l_np_{n-2}(x),\end{equation}  where $C_n$ and $\l_n$ do not depend on $x$, $p_{-1}\equiv0,$ $p_0\equiv1$ and $\l_n>0$,  and that the zeros of $p_n$ and $p_{n-1}$ interlace, i.e.,
\begin{equation}\nonumber
  x_{n,1}<x_{n-1,1}<x_{n,2}<\cdots <x_{n,n-1}<x_{n-1,n-1}<x_{n,n}.
\end{equation}
It is also known that, if $p_n$ and $p_{n-2}$ do not have a common zero, the $n-1$ zeros of $(x-C_n)p_{n-2}(x)$ interlace with the $n$ zeros of $p_n$  \cite[Theorem 3]{Beardon}, i.e., $x_{n,1}<C_{n}<x_{n,n}$ and the point $C_n$ is a natural inner bound for the extreme zeros of $p_n$.

By iterating \eqref{TTRR1}, Beardon  \cite[Theorem 4]{Beardon} generalised this result and obtained recurrence equations involving polynomials $p_{n-m}$, $p_n$ and $p_{n-1}$:
\begin{equation} \l_n\l_{n-1}\dots \l_{n-m+2}p_{n-m}(x)=S_{m-1}(x)p_{n-1}(x)+S_{m-2}(x)\label{Beardon1}p_{n}(x),m\ge2,\end{equation}
where the polynomials $S_m$ are of degree $m$ and the $n-1$ zeros of $S_{m-1}p_{n-m}$ interlace with the $n$ zeros of $p_n$, a phenomenon that will be referred to as ``completed Stieltjes interlacing''.

In this paper we re-phrase \cite[Theorem 4]{Beardon} and we introduce the associated polynomials of the orthogonal sequence $\{p_n\}_{n=0}^\infty$, as well as the three-term recurrence  satisfied by them  (cf. \cite{chiharabook,DeB,Vin}). Then we show how Christoffel's formula \cite[Theorem 2.5]{Szego_1975} can be used to obtain mixed three-term recurrence equations, similar to \eqref{Beardon1}, involving polynomials $p_n$, $p_{n-1}$, and $g_{n-m,k},m\in\{2,3,\ldots\}, k $ an integer, where the polynomial $g_{n-m,k}$ belongs to a related sequence, orthogonal with respect to the weight $c_k(x)w(x)>0$ and $c_k$ is a polynomial of degree $k$.  Restrictions on the integer $k$ are necessary, in order for the zeros of $g_{n-m,k}$ and $p_n$ to interlace in the Stieltjes sense. Each mixed recurrence equation provides us with a polynomial $G_{m-1},m\in\{2,3,\ldots,n-1\}$, whose zeros can be used as inner bounds for the extreme zeros of $p_n$. The bounds obtained in this way are more accurate than the inner bounds obtained using mixed recurrence equations with $m=2$, as was done for the extreme zeros of the Jacobi, Laguerre and Gegenbauer polynomials in \cite{Driver_Jordaan_2012}, Meixner and Krawtchouk polynomials in \cite{Jooste_Jordaan} and Hahn polynomials in \cite{JNK_2017}. In our applications the polynomials $g_{n-m,k},m\in\{2,3,\ldots,n-1\},$ are typically obtained from the polynomials in the orthogonal sequence $\{p_n\}_{n=0}^{\infty}$ by making appropriate parameter shifts of (in total) $k$ units.

Levit \cite{Levit} was the first to study the separation of the zeros of different sequences of Hahn polynomials and interlacing results  for Jacobi polynomials \cite{Askey,DJM},  Krawtchouk polynomials \cite{CS,JT} and Meixner and Meixner-Pollaczek polynomials \cite{JT} followed. Interlacing results for the zeros of different sequences of $q$-orthogonal sequences with shifted parameters are given in \cite{Gochhayat_et_al_2016, Jordaan_Tookos_2010, Moak, KoepfJT_2017}. Completed Stieltjes interlacing of zeros of different orthogonal sequences was done for the Gegenbauer  \cite{Geg}, Laguerre  \cite{Laguerre} and  Jacobi polynomials \cite{DJJ} and apart from the papers cited in the previous paragraph, inner bounds for the extreme zeros of  Gegenbauer, Laguerre  and  Jacobi polynomials were also given in \cite{Area_et_al_2012, Bottema_1931, Gupta_Muldoon_2007, Krasikov_2006, Neumann_1921, Szego_1975}; bounds for the extreme zeros of the discrete orthogonal Charlier, Meixner, Krawtchouk and Hahn polynomials in \cite{Area_et_al_2013, Krasikov_Zarkh_2009}, for the extreme zeros of the $q$-Jacobi and $q$-Laguerre polynomials in \cite{Gupta_Muldoon_2007} and for the little $q$-Jacobi polynomials in \cite{Gochhayat_et_al_2016}. Lower bounds for $x_{n,1}$ and upper bounds for $x_{n,n}$ can be found in the case of classical continuous and discrete orthogonal polynomials in \cite{Area_et_al_2004, Area_et_al_2013, Dimitrov_Nikolov_2010, Dimitrov_Rafaeli_2009,   Ismail_Li_1992, Krasikov_2002, Krasikov_2006, Krasikov_Zarkh_2009,  Szego_1975} and in \cite{Kras2005}, non-asymptotic bounds on the extreme zeros of (symmetric) orthogonal polynomials are given in terms of the coefficients of their three-term recurrence equations.

In general it is time-consuming to find the zeros of an expanded polynomial of high degree using a computer program such as Maple, since one has to work with high precision. In such cases it is helpful to have a formula for a bound.

Our main results are stated in section 2, where we also describe the algorithmic approach  used to obtain the mixed recurrence equations necessary to prove Stieltjes interlacing between the zeros of polynomials from different orthogonal sequences.  In the last three sections we display mixed recurrence equations that provide formulas for inner bounds for the extreme zeros of some specific polynomial systems, merely to show the accuracy of the bounds provided by these recurrence equations. In section 3 we show bounds in the generalized Laguerre case. In section 4 we show the quality of bounds obtained for the extreme zeros of the little $q$-Jacobi and the alternative $q$ Charlier (or $q$-Bessel) polynomials, where the equation providing the bounds concerns polynomials  $p_n(x;\a,\b)$, $p_{n-1}(x;\a,\b)$ and either $p_{n-m}(x;\a+k,\b)$ or $p_{n-m}(x;\a,\b+k)$, i.e., only one parameter is shifted.  Since we observe that the best possible bounds are obtained when one parameter is shifted optimally,  polynomial systems depending on 2 or 3 parameters also belong to this group. In section 5, we provide bounds for the extreme zeros of the  Stieltjes-Wigert and the discrete $q$-Hermite II  polynomials, as examples of  polynomial systems where no parameters are available to be shifted.

On the webpage  \url{http://www.mathematik.uni-kassel.de/~koepf/Publikationen}, we provide optimal (upper and/or lower) inner bounds for the extreme zeros of following polynomial systems:
\begin{itemize}
\item[(i)] The Gegenbauer, Laguerre, Jacobi, Bessel, Meixner and Hahn polynomials, as well as the  big $q$-Jacobi, $q$-Hahn, $q$-Meixner, little $q$-Laguerre (or Wall), little $q$-Jacobi, $q$-Krawtchouk,  affine  $q$-Krawtchouk, $q$-Laguerre, alternative $q$-Charlier (or $q$-Bessel) and $q$-Charlier polynomials. In these cases the bounds were obtained from shifting one parameter optimally. By shifting $\a$ to $\a q^{-6}$ and  $\gamma$ to $\gamma q^{-6}$, respectively, we obtain  a lower bound for the largest zero of the $q$-Charlier polynomials and  an upper bound for the smallest zero of the $q$-Meixner polynomial;
\item [(ii)] The Hermite, Krawtchouk and Charlier polynomials, as well as the quantum $q$-Krawtchouk, Stieltjes Wigert, Al-Salam Carlitz I and II and discrete $q$-Hermite I and II polynomials, where we use equations similar to \eqref{Beardon} to obtain the bounds. In this way we also find a lower bound for $x_{n,n}$ in the $q$-Meixner case and an upper bound for the smallest zero of the $q$-Charlier polynomials. 
\end{itemize}

\section{Stieltjes interlacing of zeros of classical orthogonal sequences}
\begin{lemma}\label{BeardonTh}(cf. \ \cite[Theorem 4]{Beardon})
Suppose $\{p_n\}_{n=0}^{\infty}$ is a sequence of polynomials, satisfying \eqref{TTRR1}. Then, given $n$, there exists a sequence of real orthogonal polynomials $S_m^{(n)}(x), m\in\{0,1,2,\dots,n\},$ of exact degree $m$,
satisfying the three-term recurrence equation
\begin{equation}S_{m}^{(n)}(x)=\label{TTRR2}(x-C_{n-(m-1)})S_{m-1}^{(n)}(x)-\l_{n-(m-2)}S_{m-2}^{(n)}(x),m\in\{1,2,\dots,n\},\end{equation}
 with $S_0^{(n)}(x)=1$, $S_{-1}^{(n)}(x)=0$, $S_1^{(n)}(x)=x-C_n$,
such that, for $m\ge 2$,
\begin{equation}\l_{n}\l_{n-1}\ldots \l_{n-m+2}p_{n-m}(x)=S_{m-1}^{(n)}(x)p_{n-1}(x)-S_{m-2}^{(n-1)}(x)\label{Beardon}p_{n}(x).\end{equation}
\end{lemma}

\begin{proof}
For each $n$ we have
$$\l_n \left[ {\begin{array}{cc}
   p_{n-2}(x) \\
   p_{n-1}(x)\\
  \end{array} } \right]
= \left[ {\begin{array}{cc}
   x-C_n & -1 \\
   \l_n  & 0 \\
  \end{array} } \right]
  \left[ {\begin{array}{cc}
   p_{n-1}(x)  \\
   p_{n}(x)  \\
  \end{array} } \right].$$
Thus
\begin{eqnarray}&&\nonumber\l_n \l_{n-1}\ldots\l_{n-(m-2)}\left[ {\begin{array}{cc}
   p_{n-m}(x) \\
   p_{n-(m-1)}(x)\\
  \end{array} } \right]\\
&=&\label{1}\left[ {\begin{array}{cc}
   x-C_{n-(m-2)} & -1 \\
   \l_{n-(m-2)}  & 0 \\
  \end{array} } \right]\dots
\left[ {\begin{array}{cc}
   x-C_{n-1} & -1 \\
   \l_{n-1}  & 0 \\
  \end{array} } \right]
\left[ {\begin{array}{cc}
   x-C_n & -1 \\
   \l_n  & 0 \\
  \end{array} } \right]
  \left[ {\begin{array}{cc}
   p_{n-1}(x)  \\
   p_{n}(x)  \\
  \end{array} } \right]\\
&=&\nonumber\prod_{j=0}^{m-2} \left[ {\begin{array}{cc}
   x-C_{n-j} & -1 \\
   \l_{n-j}  & 0 \\
  \end{array} } \right]
  \left[ {\begin{array}{cc}
   p_{n-1}(x)  \\
   p_{n}(x)  \\
  \end{array} } \right]. \end{eqnarray}
  Suppose
$$\left[ {\begin{array}{cc}
   A_m(x) & B_m(x)\\
   C_m(x)  & D_m(x)\\
  \end{array} } \right]=\prod_{j=0}^{m-2} \left[ {\begin{array}{cc}
   x-C_{n-j} & -1 \\
   \l_{n-j}  & 0 \\
  \end{array} } \right]$$
  Then
 $$\left[ {\begin{array}{cc}
   A_2(x) & B_2(x)\\
   C_2(x)  & D_2(x)\\
  \end{array} } \right]= \left[ {\begin{array}{cc}
   x-C_{n} & -1 \\
   \l_{n}  & 0 \\
  \end{array} } \right]$$
  and
 $$\left[ {\begin{array}{cc}
   A_3(x) & B_3(x)\\
   C_3(x)  & D_3(x)\\
  \end{array} } \right]= \left[ {\begin{array}{cc}
   x-C_{n-1} & -1 \\
   \l_{n-1}  & 0 \\
  \end{array} } \right]\left[ {\begin{array}{cc}
   A_2(x) & B_2(x)\\
   C_2(x)  & D_2(x)\\
  \end{array} } \right],$$
  from which it follows that
  $A_2(x)=x-C_n$ and
  $$A_3(x)= (x-C_{n-1})A_2(x)-C_2(x)=(x-C_{n-1})A_2(x)-\l_n.$$

In general,
 $$\left[ {\begin{array}{cc}
   A_{m}(x) & B_{m}(x)\\
   C_{m}(x)  & D_{m}(x)\\
  \end{array} } \right]= \left[ {\begin{array}{cc}
   x-C_{n-(m-2)} & -1 \\
   \l_{n-(m-2)}  & 0 \\
  \end{array} } \right]\left[ {\begin{array}{cc}
   A_{m-1}(x) & B_{m-1}(x)\\
   C_{m-1}(x)  & D_{m-1}(x)\\
  \end{array} } \right],$$
  from which we obtain
  $A_{m}(x)= (x-C_{n-(m-2)})A_{m-1}(x)-C_{m-1}(x)$ and
  $ C_{m}(x)= \l_{n-(m-2)} A_{m-1}(x),$
  i.e.,
  $$A_{m}(x)= (x-C_{n-(m-2)})A_{m-1}(x)-\l_{n-(m-3)} A_{m-2}(x).$$
  Since $A_m$ is a polynomial of degree $m-1$, we introduce the sequence of polynomials $\{S_m^{(n)}\}_{m=0}^n,$ such that $A_{m+1}(x)=S_{m}^{(n)}(x)$, and these polynomials satisfy  the three-term recurrence equation
  $$S_{m}^{(n)}(x)= (x-C_{n-(m-1)})S_{m-1}^{(n)}(x)-\l_{n-(m-2)} S_{m-2}^{(n)}(x),$$
   with $S_0^{(n)}(x)=1$, $S_{-1}^{(n)}(x)=0$, $S_1^{(n)}(x)=x-C_n$.

In the same way, we see that $B_m(x)$ is a polynomial of degree $m-2$, with $B_1(x)=0$ and $B_2(x)=-1$, satisfying the equation
$$B_m(x)=(x-C_{n-(m-2)})B_{m-1}(x)-\l_{n-(m-3)}B_{m-2}(x),m\ge 3,$$
and we can set $B_m(x)=-S_{m-2}^{(n-1)}(x)$.

Equation (\ref{1}) thus becomes
\begin{equation}\nonumber\l_n \l_{n-1}\ldots\l_{n-(m-2)}\left[ {\begin{array}{cc}
   p_{n-m}(x) \\
   p_{n-(m-1)}(x)\\
  \end{array} } \right]\\
=\nonumber \left[ {\begin{array}{cc}
   S_{m-1}^{(n)}(x) & -S_{m-2}^{(n-1)}(x)\\
   \l_{n-(m-2)}S_{m-2}^{(n)}(x)  & -\l_{n-(m-2)}S_{m-3}^{(n-1)}(x) \\
  \end{array} } \right]
  \left[ {\begin{array}{cc}
   p_{n-1}(x)  \\
   p_{n}(x)  \\
  \end{array} } \right] \end{equation}
and the result follows.
\end{proof}
For completeness, we state the following result, which is proved as part of  \cite[Theorem 4]{Beardon}.
\begin{corollary}\label{Cor}If $p_{n-m}(x)$ and $p_{n}(x)$ do not have any common zeros, the $n-1$ zeros of  $S_{m-1}^{(n)}(x)p_{n-m}(x)$ interlace with the $n$ zeros of $p_{n}(x).$
\end{corollary}
The polynomials $\{S_{m}^{(n)}\}_{m=0}^n,n\in\{0,1,\dots\},$  will be called  the \emph{associated} polynomials of the orthogonal sequence $\{p_n\}_{n=0}^\infty$. These polynomials are completely determined by the coefficients in \eqref{TTRR1}  and, according to Favard's Theorem \cite[Theorem 4.4]{chiharabook}, they are part of an orthogonal sequence. Furthermore, when $m=n$  in (\ref{Beardon}), it follows from Corollary \ref{Cor} that  the $n-1$ zeros of $S_{n-1}^{(n)}$ interlace, in the same way as the $n-1$ zeros of $p_{n-1}$, with the $n$  zeros of $p_{n}$ (cf. \cite[Theorem 4.1]{chiharabook}).

\begin{remark}
\begin{itemize}
\item[(i)] From (\ref{Beardon}) it is clear that, should $p_{n-m}$ and $p_n$ have any common zeros, these will also be zeros of $S_{m-1}^{(n)}$ and there can be a maximum of $m-1$ such common zeros (cf. \cite{Driver_Jordaan_2012,Gibson}).
\item[(ii)] The polynomials $S_n^{(n+1)}$ are the same as the numerator polynomials $p_n^{(1)}$ in \cite[p. 86, Definition 4.1]{chiharabook} and the associated polynomials in \cite{Vin}.
\item[(iii)] By expanding and re-writing \eqref{TTRR2}, we obtain a three-term recurrence equation involving the polynomials $S_m^{(n)}(x)$, $S_{m-1}^{(n-1)}(x), S_{m-2}^{(n-2)}(x), m\in\{1,2,\dots,n\},$  given by
\begin{equation}S_{m}^{(n)}(x)=\label{TTRR3}(x-C_{n})S_{m-1}^{(n-1)}(x)-\l_{n}S_{m-2}^{(n-2)}(x),\end{equation}
from which it is clear that the polynomials $S_n^{(n)}$ and $p_n$ are the same polynomials.

\end{itemize}
\end{remark}

Fix $k\in\mathbb{N}_{0}$  and $m\in\{2,3,\ldots,n-1\}$ and let $\{p_n\}_{n=0}^{\infty}$ be a sequence of orthonormal polynomials associated with the weight function $w(x)>0$ on the (finite or infinite) interval $(a,b)$. Then, from Christoffel's formula  \cite[Theorem 2.5]{Szego_1975}, the sequence of polynomials $\{g_{n,k}\}_{n=0}^{\infty}$,
orthogonal with respect to $c_k(x)w(x)>0$ on $(a,b)$, where $c_{k}(x)$ is a polynomial of degree $k$, satisfies
\begin{eqnarray}
c_k(x)g_{n-m,k}(x)&=&   \begin{array}{|cccc|}\label{det}
                    p_{n-m}(x) & p_{n-m+1}(x) & \dots & p_{n-m+k}(x) \\
                     p_{n-m}(x_1)& p_{n-m+1}(x_1) & \dots & p_{n-m+k}(x_1)\\
                      \dots & \dots & \dots & \dots  \\
                      p_{n-m}(x_k) &  p_{n-m+1}(x_k) & \dots & p_{n-m+k}(x_k)\\
                      \end{array}\\
               &=& \nonumber \sum_{j=0}^{k} U_j p_{n-m+j}(x),
               \end{eqnarray}
where $x_j,j\in\{1,2,\dots,k\}$ are the zeros of $c_k$ and the coefficients $U_j,j\in\{0,1,\dots,k\}$ are determinants. We note that, in the case of a zero $x_j$ of multiplicity $s>1$, the corresponding rows of the determinant (\ref{det}) are replaced by the derivatives of order $0,1,\dots,s-1$ of the polynomials $p_{n-m}(x), p_{n-m+1}(x),\ldots,$ $p_{n-m+k}(x)$ at $x=x_j$.

Furthermore, using Lemma \ref{BeardonTh},  the $(m-1)$ polynomials $p_{n-m+j},j\in\{0,1,\dots,m-2\},$ can be written in terms of the polynomials $p_n$ and $p_{n-1}$, by using
$$\lambda_n \lambda_{n-1}\ldots \lambda_{n-(m-j-2)}p_{n-(m-j)}(x)=S_{m-(j+1)}^{(n)}(x)p_{n-1}(x)-S_{m-(j+2)}^{(n-1)}(x)p_{n}(x).$$
In order to write the rest of the polynomials $p_{n-m+j},j\in\{m+1,m+2,\dots,k\},$ in terms of the polynomials $p_n$ and $p_{n-1}$, we use the following iterations of the three-term recurrence equation (\ref{TTRR1})
\begin{equation}p_{n+m}(x)=\nonumber S_m^{(n+m)}(x)p_{n}(x)-\l_{n+1}S_{m-1}^{(n+m)}(x)p_{n-1}(x),\end{equation}
where the polynomials $S_{m}^{(n+m)}(x), m\in\{1,2,3,\dots\},$  satisfy (\ref{TTRR3}), with $n$ replaced by $n+m$.
Or, when we replace $m$ by $-m+j$, we get
\begin{equation}p_{n-m+j}(x)=\nonumber S_{-m+j}^{(n-m+j)}(x)p_{n}(x)-\l_{n+1}S_{-m+j-1}^{(n-m+j)}(x)p_{n-1}(x)\end{equation}
and (\ref{det}) can be written as:
\begin{eqnarray}
c_k(x)g_{n-m,k}(x)&=&\nonumber\Bigg(\sum_{j=0}^{m-2}\frac{U_j~S_{m-(j+2)}^{(n-1)}(x)}{\lambda_n \lambda_{n-1}\ldots \lambda_{n-(m-j-2)}}+U_{m}+\sum_{j=m+1}^{k}U_j~S_{j-m}^{(n-m+j)}(x)\Bigg)p_{n}(x)\\
&+&\nonumber \Bigg(\sum_{j=0}^{m-2}\frac{U_j~S_{m-(j+1)}^{(n)}(x)}{\lambda_n \lambda_{n-1}\ldots \lambda_{n-(m-j-2)}}+U_{m-1}-\lambda_{n+1}\sum_{j=m+1}^{k}U_j~S_{j-m-1}^{(n-m+j)}(x)\Bigg)p_{n-1}(x)\\
&=&\nonumber R(x)p_{n}(x) +  G(x)p_{n-1}(x)
\end{eqnarray}
Christoffel's formula thus provides us with a mixed three-term recurrence equation, involving polynomials $g_{n-m,k}(x),p_n$ and $p_{n-1}$ and the coefficients in the mixed three-term recurrence equation are determined by the coefficients in \eqref{TTRR1} and the determinants $U_j, j\in\{0,1,\dots,k\}$. It is clear that the coefficient $G(x)$ of $p_{n-1}(x)$ is of degree $\max\{m-1,k-m-1\}$ and the coefficient $R(x)$ of $p_{n}(x)$ is of degree $\max\{m-2,k-m\}$.

In order for the zeros of $G(x)g_{n-m,k}(x)$ to fully interlace with the zeros of $p_n(x)$, the polynomial $G(x)g_{n-m,k}(x)$ must be of exact degree $n-1$, which means that $G(x)$ must be of exact degree $m-1$, i.e., $\max\{m-1,k-m-1\}=m-1,$ i.e., $k-m-1\le m-1$, from which we can deduce  $k\le 2m$.

We have thus proved the following result.

\begin{theorem}\label{main2} Let $\{p_n\}_{n=0}^{\infty}$ be a sequence of polynomials orthogonal on the (finite or infinite) interval $(a,b)$ with respect to the weight function $w(x)>0.$ Let $k\in\mathbb{N}_{0}$  be fixed. Suppose $\{g_{n,k}\}_{n=0}^{\infty}$ is a sequence of polynomials
orthogonal with respect to $c_k(x)w(x)>0$ on $(a,b)$, where $c_{k}(x)$ is a polynomial of degree $k$. Then, for $m\in\{2,3,\ldots,n-1\}$ and $k\in\{0,1,\dots,2m\},$ there exist  on $(a,b)$ a polynomial $G_{m}(x)$ of degree $m$, and a polynomial $a_{k-m}$  that is of degree $m-2$ when $k-m\in\{-m,-m+1,\ldots,m-2\}$  and of degree $k-m$ whenever $k-m\in\{m-1,m,m+1,\ldots\}$, such that
\begin{equation}\label{general4}A_nc_{k}(x) g_{n-m,k}(x)=a_{k-m}(x) p_{n}(x)-G_{m-1}(x) p_{n-1}(x), ~n\in\{3,4,\dots\}.\end{equation}
Furthermore,
\begin{itemize}
\item[(i)] if $g_{n-m,k}$ and $p_n$ are co-prime, the $n-1$ real, simple zeros of $G_{m-1}(x)g_{n-m,k}$ interlace with the zeros of $p_n$, the smallest zero of $G_{m-1}$ is an upper bound for the smallest zero of $p_n$, and the largest zero of $G_{m-1}$ is  a lower bound for the largest zero of $p_n$;
\item[(ii)] if $g_{n-m,k}$ and $p_n$ are not co-prime and have $r$ common zeros counting multiplicity, then\newline
\noindent a) $r\le \min\{m,n-m-1\}$; \newline
\noindent b) these  $r$ common zeros are simple zeros of $G_{m-1}$;\newline
\noindent c) no two successive zeros of $p_n$, nor its largest or smallest zero can be a zero of $G_{m-1}$;\newline
\noindent d) the $n-2r-1$ zeros of $G_{m-1}g_{n-m,k}(x),$  none of which is a zero of $p_n$, together with the $r$ common zeros of $g_{n-m,k}$ and $p_n$, interlace with the $n-r$ non-common zeros of $p_n$; \newline
\noindent e) the smallest zero of $G_{m-1}$ is an upper bound for the smallest zero of $p_n$, and the largest zero of $G_{m-1}$ is  a lower bound for the largest zero of $p_n.$
\end{itemize}
\end{theorem}
For the proofs of (i) and (ii) above, we refer the reader to \cite[Theorem 2.1,  Corollary 2.2]{Driver_Jordaan_2012}.
\begin{remark} In  \cite{Jooste_Jordaan} the authors use equations like \eqref{general4}, with $m=2$, to obtain inner bounds for the Meixner polynomials. In the proof of \cite[Theorem 2.1]{Jooste_Jordaan}, the assumption is made that if $q$ is any polynomial, such that
 $$\int_{a}^{b}q(x)p_{n}(x)w(x)dx=0,n\in\{3,4,\dots\},$$
then $\mbox{deg}(q(x))\le n-1$, which is not always true. \end{remark}

In the following sections we will assume that the polynomials $g_{n-m,k}$ and $p_n$ involved in \eqref{general4} do not have any common zeros, in order to obtain inner bounds for the extreme zeros of $p_n$. Due to their complexity, we use computer algebra to find equations  like \eqref{general4} and, following the approach in \cite{Chen, Koepf2014}, we write a procedure that applies the Gosper  \cite{Gosper_1978, Koepf2014} and Zeilberger \cite{Koepf1995, Koepf2014, Petkovsek_et_al_1996, DZ} algorithms.

Gosper's algorithm deals with finding, for $j$ an integer, an anti-difference  $s_j,$  for given $a_j$, i.e., a sequence $s_j$ for which $a_j=\Delta s_j=s_{j+1}-s_j$, in the particular case that $s_j$ is a hypergeometric term, i.e.,
 \[\displaystyle\frac{s_{j+1}}{s_j}\in \mathbb{Q}(j).\]
Given the hypergeometric term $F(n,j,\a)$  with respect to variables $n,j$ and $\a$, Zeilberger's algorithm  provides a  recurrence equation (with polynomial coefficients) for \[ s_n=\sum_{j=-\infty}^\infty F(n,j,\a).\]
Suppose we have the particular case where
   \[ a_j=F(n-m,j,\alpha+k)+\sum_{i=0}^J\sigma_i(n)F(n-i,j,\alpha),\]
    with undetermined variables $\sigma_i(n)$, $i\in\{0,1,\ldots,J\},$ and
     \[s_n:=p_n(x;\a)=\sum_{j=-\infty}^\infty F(n,j,\alpha).\]
      We apply Gosper's algorithm to $a_j$ and, if successful, an anti-difference $g(n,j)$ of $a_j$ is found, i.e., $a_j =g(n,j+1)-g(n,j)$ and at the same time $\sigma_i(n)$, $i\in\{0,1,\ldots,J\}$. By summation, we obtain
\begin{align*}
 &\sum_{j=-\infty}^\infty a_j=\sum_{j=-\infty}^\infty \left( F(n-m,j,\alpha+k)+\sum_{i=0}^J\sigma_i(n)F(n-i,j,\alpha)\right)  \\
 & =p_{n-m}(x;\a+k)+\sum_{i=0}^J\sigma_i(n)p_{n-i}(x;\a)= \sum_{j=-\infty}^\infty \left(g(n,j+1)-g(n,j) \right)=0,
    \end{align*}
since the last sum is telescoping. It follows that $p_{n}(x;\a)$ satisfies a recurrence equation of type \eqref{general4}. We refer the reader to \cite[Chapters 5--7]{Koepf2014} and references there-in for more details about the algorithms of Gosper and Zeilberger implemented in the Maple \texttt{hsum} package. The $q$-analogues of Gosper's and Zeilberger's algorithms  are implemented in the Maple \texttt{qsum} package \cite{Koepf2014}. We apply an adaption of the \texttt{sumdiffeq} \cite[p.~210]{Koepf2014} and \texttt{qsumdiffeq} \cite[p.~219]{Koepf2014} procedure of the \texttt{hsum} and \texttt{qsum} packages, respectively, in order to obtain two procedures to derive recurrence equations of type \eqref{general4} for the polynomial systems considered in the sequel. The \texttt{hsum} and the \texttt{qsum} packages can be downloaded from  \url{http://www.mathematik.uni-kassel.de/~koepf/Publikationen}, as well as the two Maple codes used to obtain our mixed three-term recurrence equations.
The first program called  \texttt{Mixedrec$(F,k,S(n),s_0,a,s_1,b,s_2)$} finds a recurrence equation of the form
\[S(n-s_0,a+s_1,b+s_2)=\sum_{j=0}^J\sigma_jS(n-j, a,b), J\in\{1,2,\ldots\},\]
where  $S(n,a,b)=\displaystyle\sum_{k=-\infty}^\infty F$, $F$ is a function of $k$, $n$, $a$ and $b$, and $s_0$, $s_1$, $s_2$ are  integers. The second one denoted by  \texttt{qMixedrec$(F,k,S(n),s_0,a,s_1,b,s_2)$} is the $q$-analogue of the first one and  finds a recurrence equation of the form
\[S(n-s_0,a q^{s_1},b q^{s_2})=\sum_{j=0}^J\sigma_jS(n-j, a,b), J\in\{1,2,\ldots\}.\]

We use this algebraic method to obtain mixed three-term recurrence equations involving polynomials  $p_n(x;\alpha,\beta)$ and $p_{n-1}(x;\alpha,\beta),$ belonging to the same sequence that is orthogonal on an interval $(a,b)$ with respect to a measure $w(x;\a,\b),$ and a polynomial from a related sequence, obtained by integer shifts of the parameters $\a$ and $\b$, namely $p_{n-m}(x;\a+s,\b+t),m\in\{2,3,\ldots,n-1\}$, which is orthogonal with respect  to
\begin{equation}\nonumber w(x;\a+s,\b+t)=c_{s+t}(x;\a,\b)w(x;\a,\b)>0\end{equation}
 on $(a,b)$, where $c_{k}(x;\a,\b)$ is a polynomial of degree $k$ in $x$. If the sequence is $q$-orthogonal with respect to the weight  $w(x;\a,\b)$, the equations involve the polynomials  $p_n(x;\alpha,\beta)$ and $p_{n-1}(x;\alpha,\beta),$ and $p_{n-m}(x;\a q^s,\b q^t),m\in\{2,3,\ldots,n-1\}$ and the latter polynomial  is orthogonal with respect to $$w(x;\a q^s,\b q^t)=c_{s+t}(x;\a,\b)w(x;\a,\b)>0$$ on $(a,b)$. From Theorem  \ref{main2} we know that for positive integers $s$ and $t$ such that $s+t\in\{0,1,\ldots,2m\}$,
the polynomial coefficient of $p_{n-1}(x;\a,\b)$ in this mixed recurrence equation, is of exact degree $m-1$ and will be denoted by $G_{m-1}^{(\a+s,\ \b+t)}(x).$


\section{Bounds for the extreme zeros of the classical orthogonal polynomials}
The generalized Laguerre polynomials
\[L_n^{(\alpha)}(x)={(\alpha+1)_n\over n!}\hypergeom{1}{1}{-n}{\alpha+1}{x} ,~\alpha>-1\] are orthogonal with respect to $w(x;\a)=x^{\a}e^{-x}$ on $(0,\infty).$
Since optimal bounds are obtained (cf.\ \cite {Driver_Jordaan_2012}) when we shift $\a$ optimally, we use  mixed three-term recurrence equations involving the polynomials $L_n^{(\alpha)}(x)$, $L_{n-1}^{(\alpha)}(x)$ and $L_{n-m}^{(\alpha+2m)}(x)$ for $m=3$ and $m=4$, where $L_{n-m}^{(\alpha+2m)}(x)$ is orthogonal on $(0,\infty)$ with respect to
$w(x;\a+2m)$, where
$$\frac{w(x;\a+2m)}{w(x;\a)}=x^{2m}=c_{2m}(x),$$   to illustrate the quality of our newly found bounds in the classical case.

When we take $m=3$, the mixed three-term recurrence equation
\begin{equation}\label{Lag3}
x^6L_{n-3}^{(\alpha+6)}(x)=-n a_3(x) L_{n}^{(\alpha)}(x)+(n+\alpha)G_2^{(\a+6)}(x)L_{n-1}^{(\alpha)}(x),
\end{equation}
with $$a_3(x)=\left( n-1 \right)  \left( n-2 \right)( {x}^{3}+3\,  \left( \alpha+3 \right) {x}^{2})-   \left( \alpha+2 \right)_3
\left( \alpha+4\,n-3 \right) x+ \left( \alpha+1 \right)_5$$ and
$$G_2^{(\a+6)}(x)=(\alpha+3)(3n(n+\alpha+1)+(\alpha+1)_2)x^2
-2(\alpha+2)_3(\alpha+2n+1)x+(\alpha+1)_5.$$
is clearly in the form of \eqref{general4} and from Theorem \ref{main2} we know that the $n-1$ zeros of $G_2^{(\a+6)}(x)L_{n-3}^{(\alpha+6)}(x)$ interlace with the $n$ zeros of $L_{n}^{(\alpha)}(x)$. The zeros of $G_2^{(\a+6)}(x)$ are
\begin{eqnarray}\label{bound_lag}
B_{3}^{(\a+6)}&=&\frac{(\alpha+2)(\alpha+4)(\alpha+2n+1)}{3n(n+\alpha+1)+ \left( \alpha+1 \right)_2 }\\
 &\pm &\nonumber\frac{\sqrt{(\alpha+2)(\alpha+4)\Big(\left( {\alpha}^{2}+6\,\alpha+17 \right)( {n}^{2}+ \alpha n +n)  - \left( \alpha+2 \right)  \left( \alpha+1 \right) ^{2}
\Big)} }{ 3n(n+\alpha+1)+ \left( \alpha+1 \right)_2 }
\end{eqnarray}
and they are inner bounds for the extreme zeros of $L_{n}^{(\alpha)}(x)$. By computing the values of these bounds, we find that the smallest value in \eqref{bound_lag} is an accurate upper bound for $x_{n,1}$, however, by substituting $m=4,k=2m=8$ in \eqref{general4}, we obtain
\begin{equation}
\nonumber x^8L_{n-4}^{(\alpha+8)}(x)=a_4(x)L_{n}^{(\alpha)}(x)+(n+\alpha)G_{3}^{(\a +8)}(x)L_{n-1}^{(\alpha)}(x),
\end{equation}
where
\begin{eqnarray}
   a_4(x)&=&\nonumber \left( n-3 \right)_4   \left( x+4\a+16 \right){x}^{3} -n\left( \alpha+3 \right)_3  \Big( 5\,{n}(2n+\alpha-4) +{\alpha}^{2}-2\,\alpha+12 \Big) {x}^{2}\\
   &+&\nonumber2\,n\left( \alpha+2 \right)_5  \left( \alpha+3\,n-2 \right) x -n\left(\alpha+1 \right)_7
\end{eqnarray}
and
\begin{eqnarray}   G_{3}^{(\a +8)}(x)&=&\label{lag4}
 - \left( \alpha+4 \right)  \left( 2\,n+\alpha+1 \right)  \left( {\alpha}^{2}+2\,\alpha\,n+2\,{n}^{2}+5\,
\alpha+2\,n+6 \right) {x}^{3}\\
&+&\nonumber\left( \alpha+3 \right)_3\left( 3(\alpha+1)_2+ 10\,n\,(n+\a+1) \right) {x}^{2}-3\, \left( \alpha+2 \right)_5 \left( 2\,n+\alpha+1 \right) x+\left( \alpha+1 \right)_7
\end{eqnarray}
 and the smallest zero of $ G_{3}^{(\a +8)}(x)$ is an even more accurate bound than the bound obtained from \eqref{Lag3}.  In Table \ref{LagB} we show examples indicating the accuracy of this value, computed numerically,  as upper bound for the lowest zero of $L_{n}^{(\alpha)}(x).$
To find a lower bound for the largest zero, we let $m=4$ and $k=0$  in \eqref{general4}, i.e., we don't consider any parameter shifts, and we obtain the recurrence equation
$$L_{n-4}^{(\alpha)}(x) =-{\frac {n \left( {x}^{2}-2 \left(\alpha+2n-4\right) x+{\alpha}^{2}+3\,\alpha\,n+
		3\,{n}^{2}-6\,\alpha-12\,n+11 \right)L_{n}^{(\alpha)}(x) }{ \left( \alpha-3+n \right)_3
	}}+{\frac {
	G_3^{(\a)}(x) L_{n-1}^{(\alpha)}(x) }{ \left( \alpha-3+n \right)_3   }}$$
with
\begin{eqnarray}
G_3^{(\a)}(x)&=&\label{Laglb}-{x}^{3}+ 3\left( \alpha+2\,n-3 \right) {x}^{2}- \left(3\,{\alpha}
^{2}+10\,\alpha\,n+10\,{n}^{2}-15\,\alpha-30\,n+18 \right) x\\
&+&\nonumber \left(
\alpha-3+2\,n \right)  \left( {\alpha}^{2}+2\,\alpha\,n+2\,{n}^{2}-3\,
\alpha-6\,n+2 \right).
\end{eqnarray}
The largest zero of $G_3^{(\a)}(x)$ is a lower bound for the largest zero of $L_{n}^{(\alpha)}(x)$. Here we use a zero of a third degree polynomial, that can easily be computed numerically, to approximate a zero of a polynomial of much higher degree.

We show the quality of these bounds in Table \ref{LagB}
and we compare our upper bound for $x_{n,1}$ with
the upper bound
\begin{equation}
\label{BoundGM} B_{GM}={\frac { \left( \alpha+1 \right)  \left( \alpha+2 \right)  \left(\alpha+4 \right)  \left( 2\,n+\alpha+1 \right) }{ \left( \alpha+1 \right) ^{2} \left( \alpha+2 \right) + \left( 5\,\alpha+11 \right) n \left( n+\alpha+1 \right) }}
\end{equation}
 obtained in \cite[Equation 2.11]{Gupta_Muldoon_2007}. We also provide the lower bound for the largest zero, as  given in \cite{Krasikov_2006} for values of  $n\ge 30$. This bound is more precise than the bound obtained from (\ref{Laglb}), however, the recurrence equation involving polynomials $L_{n}^{(\alpha)}$, $L_{n-1}^{(\alpha)}$ and $L_{n-7}^{(\alpha)}$, which is too big to display, provides us with a polynomial $G_6^{(\a)}(x)$ and the largest zero of this polynomial, that can be computed numerically, is a more accurate bound than the bound in \cite{Krasikov_2006}. The equation can be found in the accompanying Maple file.
 \begin{table}[h!]\centering \caption{Examples to show the quality of the inner bounds for the extreme zeros of $L_{n}^{(\alpha)}(x)$ for different values of $n$ and  $\alpha$.}\label{LagB}
	{\small{\begin{tabular}{|c|c|c|c|c|c|c|c|c|}
		\hline
		$n$&$\alpha$&$x_{n,1}$&\small{Bound from \eqref{lag4}}&$B_{GM}$ \eqref{BoundGM}&\footnotesize{Bound from (\ref{Laglb})}& \footnotesize{Bound in \cite{Krasikov_2006}} &\small{$x_{n,n}$}\\\hline
10&-0.5&0.06019206315&0.06019206332&0.060269&28.469&n/a&29.025\\\hline
		10&10.0&3.50805&3.51556&4.0168&44.945&n/a&46.365\\\hline
		100&10.0&0.49692&0.49863&0.5746&353.225&387.960&394.294\\\hline
		100&-0.5&0.00615313229&0.00615313230&0.0061611&335.472&367.816&374.006\\\hline
	\end{tabular}}}
\end{table}		
Similar results can be obtained for the Jacobi, Gegenbauer, Hermite, Bessel, Meixner and Hahn polynomials and the recurrence equations providing the bounds, as well as the bounds, are also available in the Maple file.
\begin{remark}\begin{itemize}
\item [(i)] The Charlier polynomials
 \[C_n(x;a)=\hypergeom{2}{0}{-n,-x}{-}{-\frac 1a}, a>0,\] are orthogonal on $(0,\infty)$ with respect to $w(x;a)=\frac{a^x}{x!}$ and when we shift the parameter $a$ by $k$ units, we obtain  the polynomial $C_n(x;a+k)$, which is orthogonal with respect to $w(x;a+k)$ on  $(0,\infty)$ and since $\frac{w(x;a+k)}{w(x;a)}=\frac{(a+k)^x}{a^x}$ is not a polynomial of degree $k$ in $x$, these polynomials do not satisfy the conditions of Theorem \ref{main2} and we use Lemma \ref{BeardonTh} and Corollary \ref{Cor} to obtain bounds for the extreme zeros of these polynomials.
 \item [(ii)] We also use Lemma \ref{BeardonTh} in the case of the Krawtchouk polynomials
 \[K_n(x;p,N)=\hypergeom{2}{1}{-n,-x}{-N}{\frac 1p}, n\in\{0,1,\dots,N\},N\in\nn,\] that are orthogonal for $0<p<1$, since integer shifts of the parameter $p$ will result in the parameter moving outside the interval where orthogonality is guaranteed.
 \end{itemize}
 \end{remark}
\section{Bounds for the extreme zeros of $q$-orthogonal polynomials obtained by shifting one parameter}
In this section we illustrate the method by obtaining new bounds for the extreme zeros of the  little $q$-Jacobi and the alternative $q$-Charlier (or $q$-Bessel) polynomials.

\subsection{The little $q$-Jacobi polynomials}
The little $q$-Jacobi polynomials (cf.\ \cite[Section 14.12]{KLS})
\[p_n(q^x;\alpha,\beta|q)=\qhypergeom{2}{1}{q^{-n},\alpha\beta q^{n+1}}{\alpha q}{q;q q^x},0<\alpha q<1,\, \beta q<1\] are discrete orthogonal with respect to the weight $\displaystyle w(x;\a,\b)=\frac{(\b q;q)_x}{(q;q)_x}(\a q)^x$ on the interval $(0,1).$ The polynomial $p_{n}(q^x;\alpha q^4,\beta|q)$ is orthogonal on $(0,1)$ with respect to $(q^x)^4 w(x;\a,\b).$ By replacing $q^x$ with $x$, we obtain
\begin{eqnarray}\nonumber
&&\frac{\alpha\, ( {q}^{n}-q) (\a \b {q}^{2\,n}-1) (\b {q}^{n}-q )(
 \alpha\,\beta\,{q}^{n+1};q)_2 ( \alpha\,{q}^{n+1};q)_2}{( \alpha\,q;q)_4} {x}^{4}p_{n-2}(x;\alpha q^4,\beta|q)  \\
&=& a_2(x)p_{n}(x;\alpha ,\beta|q)+{q}^{3\,n-5}\left(1 -\alpha{q}^{2} \right) G^{(\alpha q^4,\ \beta)}(x)p_{n-1}(x;\alpha ,\beta|q),
\end{eqnarray}
 with
 \begin{equation} \nonumber
a_2(x)=\alpha\,( {q}^{n}-q)(\b {q}^{n}-q)
  \left(  \left(\a \b {q}^{2\,n}-1 \right) {q}^{2\,n-3}x+
  \left( q+1 \right)  \left(\a {q}^{2}-1 \right) {q}^{3\,n-5}
  \right)x+(\a q;q)_3 {q}^{4\,n-6}
 \end{equation}
and $G^{(\alpha q^4,~\beta)}(x)$ is a linear function with zero
\begin{equation}
B_2^{(\alpha q^4,~\beta)}=\label{BndqJ}{\frac { \left( \alpha\,{q}^{3}-1 \right)  \left( \alpha\,q-1 \right)
{q}^{n-1}}{(\a\,\beta\,{q}^{2\,n+1}+1)({\alpha}{q}^{2}+1)-\a\,{q}^{n+1}(\b+1)(q+1)}}.
\end{equation}

In Table \ref{LagBd} we show the quality of  $B_2^{(\alpha q^4,\ \beta)}$ as bound and we compare it to bounds given in \cite[Equation 4.3]{Gupta_Muldoon_2007} when $n=10,30$ and  in \cite[Equation 4.2]{Gupta_Muldoon_2007} for $n=100$ (see Remark \ref{rem}).  We also observe that the value of $q$ is an upper bound for the zero $x_{n,n-1}$. The accuracy of this upper bound decreases for $q$ in the vicinity of $1$. Furthermore, we observe that $\frac{x_{n,j+1}}{x_{n,j}}\approx \frac 1q$ for $j\in\{1,2,\ldots,n-1\}$ and again this is less accurate when $q \rightarrow 1$.
\begin{table}[h]\centering \caption{Examples to show the quality of the inner bounds for the extreme zeros of $p_{n}{(x;0.5,\b|q)}$ for different values of $n,\b$ and $q$.}\label{LagBd}
\begin{tabular}{|c|c|c|c|c|c|c|c|}
		\hline
		$n$&$\b$&$q$&$x_{n,1}$&$B_2^{(\alpha q^4,\ \beta)}$ in \eqref{BndqJ}&Bound in \cite{Gupta_Muldoon_2007}&$x_{n,n-1}$&$x_{n,n}$\\\hline
		10&1&0.6&0.005216&0.005359&0.021776 &0.600000&1.000000\\\hline
		30&-10&0.6&$1.8642\cdot10^{-7}$&$1.9497\cdot 10^{-7}$&$7.9382\cdot 10^{-7}$ &0.600000&1.000000\\\hline
		100&-10&0.9&0.0000059&0.0000073&0.0000131 &0.8999999&1.000000\\\hline
\end{tabular}
\end{table}	
\begin{remark}\label{rem}
\begin{itemize}
\item [(i)] The little $q$-Laguerre (or Wall) polynomials are obtained from the little $q$-Jacobi polynomials, by letting $\beta=0$ and the bound $B_2^{(\alpha q^4,~\beta)}$ in \eqref{BndqJ}, with $\beta=0$,  can be used as an upper bound for the smallest zero of the little $q$-Laguerre polynomial.
\item [(ii)] In \cite{Gupta_Muldoon_2007}, Gupta and Muldoon provide bounds for the smallest zeros of the  little $q$-Jacobi polynomial $p_n((1-q)x;\alpha,\beta|q )$ and the $q$-Laguerre polynomial
$L_n^{(\alpha)}((1-q)x;q),\alpha>-1.$ Using a suitable comparison we observe that, for both these systems, the upper bounds for the smallest zeros obtained by our method are more accurate than the upper bounds obtained in \cite{Gupta_Muldoon_2007}.
\end{itemize}
\end{remark}

\subsection{The alternative $q$-Charlier or $q$-Bessel polynomials}
The alternative $q$-Charlier polynomials (cf.\ \cite[Section 14.22]{KLS})
\[y_n(q^x;\alpha,q)=\qhypergeom{2}{1}{q^{-n},-\alpha q^n}{0}{q;q\,q^x}, \alpha>0,\]
are orthogonal with respect to the weight function $\displaystyle w(x;\a)=\frac{\a^x}{(q;q)_x}q^{\binom{x+1}{2}}$ on $(0,1)$ and $y_{n}(q^x;\alpha {q}^{4},q)$ are orthogonal on $(0,1)$ with respect to $(q^x)^4 w(x;\a)$. By replacing  $q^x$ with $x$,  these polynomials satisfy
 \begin{eqnarray} \nonumber
&&  \left(- \alpha\,{q}^{n},q \right)_2 \left( {q}
   	^{n}-q \right)  \left(\a {q}^{2\,n}+q \right)( \alpha{q}^{4})  {x}^{4}  y_{n-2}(x;\alpha {q}^{4},q)\\
   & =&\nonumber {q}^{2\,n}  \left(\a \left( \a{q}^{3\,n+1}-\a{q}^{2\,n+2}
+{q}^{n+2}-\,{q}^{3} \right) {x}^{2}- \a q^n\left({q}^{2}-{q}^{n+1}-{q}^{n}+{q} \right) x+{q}^{n} \right)  y_{n}(x;\alpha,q )\\
   &-&\nonumber {q}^{2\,n}\left( \left(\a{q}^{2\,n+1} -\a{q}^{n+2}-\a{q}^{n+1}-q \right)
x+{q}^{n}
 \right) y_{n-1}(x;\alpha,q ).
 \end{eqnarray}
From the  coefficient of $y_{n-1}$,  we obtain
\begin{equation}\label{Boundqc240}
   B_2^{(\alpha q^4)}=\frac {{q}^{n-1}}{1-\alpha{q}^{n}({q}^{n}+{q}+1)},
\end{equation}
    which is a  relatively sharp upper bound for the lowest zero of $y_n$, as shown in Table \ref{AQCbounds}. In the recurrence equation involving $y_{n-3}(x;\a q^6,q),y_{n-1}(x;\a ,q)$ and $y_{n}(x;\a,q),$ the coefficient of $y_{n-1}$ is the quadratic polynomial
\begin{eqnarray}\nonumber
 G_2^{(\alpha q^6)}(x)&=&{q}^{3} \left( 1+ \left( {q}^{4\,n}+({q}^{2\,n+2}-{q}^{3\,n})({q}^{2}+{q}+1) \right) {\alpha}^{2}
+ \left( {q}^{n}+{q}^{n+2}-{q}^{2\,n}+{q}^{n+1} \right) \alpha\right) {x}^{2}\\
&+& \nonumber{q}^{n+1} \left( q+1 \right)\left( \a{q}^{2\,n}-\a{q}^{n+2}-
\alpha\,{q}^{n}-1 \right) x  +{q}^{2\,n}
\end{eqnarray}
and the lowest zero of this polynomial is
\begin{equation}\label{AQC}
B_3^{(\alpha q^6)}=\frac {-b-\sqrt{b^2-4ac}}{2a},
\end{equation}
 with $c=G_2(0)$, $b=G_2'(0)$ and $a=\frac{G_2''(0)}{2},$
which is a more accurate upper bound for $x_{n,1}$, as shown in Table \ref{AQCbounds}.  In Table \ref{AQCbounds} we also illustrate the quality of the bound
\begin{equation}\label{BoundqcTTRR}
   B_n=-{\frac {{q}^{n+1} \left( {q}^{2\,n}\alpha-\alpha\,{q}^{n+1}-\alpha\,{
   			q}^{n}-{q}^{2} \right) }{ \left( {q}^{2\,n}\alpha+q \right)  \left( {q
   		}^{2\,n}\alpha+{q}^{3} \right) }},
   \end{equation}
    obtained from the three-term recurrence equation satisfied by these polynomials, which increases for large values of $n$.
\begin{table}[h]\centering \caption{Examples to show the quality of the upper bound for the lowest zero of  $y_n(x;\alpha,q)$ for different values of $n$, $\alpha$ and $q$.}\label{AQCbounds}
\begin{tabular}{|c|c|c|c|c|c|c|c|}
  \hline
  $n$&$q$&$\alpha$&$x_{n,1}$&$B_3^{(\alpha q^6)}$ in (\ref{AQC})& $B_2^{(\alpha q^4)}$ in (\ref{Boundqc240})& $B_n$ in (\ref{BoundqcTTRR})\\
  \hline
   10&0.55&0.5&0.0045925&0.0045925&0.0045964&0.004635\\\hline
 10&0.99&10&0.06207968&0.06796546&0.08444314&0.115245\\\hline
  10&0.45&100&0.00071318&0.00071318&0.00072109&0.000941\\\hline
  70&0.45&10&$1.179406\cdot10^{-24} $ &$1.179406\cdot10^{-24}$&$1.179406\cdot10^{-24}$&$1.179406\cdot10^{-24}$\\\hline
  70&0.8&100&$2.05671\cdot10^{-7}$&$2.05671\cdot10^{-7}$&$2.05682\cdot10^{-7}$&$2.05783\cdot10^{-7}$\\\hline
  \end{tabular}
  \end{table}	

The value of $q$ is an upper bound for the zero $x_{n,n-1}$ of $y_n(x;\alpha,q)$, as shown in Table \ref{boundq} below and the accuracy of this upper bound decreases for $q$ in the vicinity of $1$. Furthermore, $\frac{x_{n,j+1}}{x_{n,j}}\approx \frac 1q$ for $j\in\{1,2,\ldots,n-1\}$ and again this is less accurate for the values of $q$ in the vicinity of $1$.
\begin{table}[h]\centering \caption{Examples to show the quality of $q$ as upper bound for $x_{n,n-1}$ for different values of $n$, $\alpha$ and $q$.}\label{boundq}
\begin{tabular}{|c|c|c|c|c|}
   \hline
   $n$&$q$&$\alpha$&$x_{n,n-1}$&$x_{n,n}$\\
   \hline
10&0.55&0.50&0.54999999&1.00000000\\\hline
   10&0.99&0.10&0.9735504&0.9935188\\\hline
   10&0.10&1000&0.099999999&1.00000000\\\hline
   70&0.70&1000&0.700000000&0.99999999\\\hline
   70&0.70&0.10&0.700000000&0.99999999\\\hline
\end{tabular}
   \end{table}		

\begin{remark}
Finding inner bounds by using Theorem \ref{main2} is not possible for the following polynomial systems, since
\begin{itemize}
\item [(i)] the quantum $q$-Krawtchouk polynomials $K_n^{(qtm)}(q^{-x};p,N;q)$ are orthogonal on $[0,N]$ for $p>q^{-N}$ and shifting $p$ causes a change in the interval of orthogonality;
 \item [(ii)] the Al-Salam Carlitz I polynomials $U_n^{(\a)}(x;q)$ are orthogonal for $\a<0$ on $(\a,1)$ and shifting $\a$ to $\a q^k$ will result in a change in the interval of orthogonality;
 \item [(iii)] the Al-Salam Carlitz II polynomials $V_n^{(\a)}(q^{-x};q)$ are orthogonal for $0<\a q<1$  on $(0,1)$ with respect to
   $$  w(x;\a)=\frac{q^{x^2}\a^x}{(q;q)_k(\a q;q)_k} \text{~~and~~}\frac{w(x;\a q^{-k})}{w(x;\a)}=\frac{\left(\frac {1}{\alpha} q^{-x};q\right)_k}{\left(\frac {1}{\alpha} ;q\right)_k}=c_k(q^{-x};\a)$$ is a polynomial of degree $k$ in the variable $q^{-x}$. However, when we substitute $\a$ with $\a q^{-k},k\in\{1,2,\dots\}$, the condition $0<\a q<1$ is not satisfied.
 \end{itemize}
 We find bounds for the zeros of these systems by using Corollary \ref{Cor}.
 \end{remark}

\section{Bounds for the extreme zeros of $q$-orthogonal polynomials obtained without any parameter shifts}
In this section we obtain inner bounds for the extreme zeros of the Stieltjes-Wigert and discrete $q$-Hermite II polynomials, where no parameter shifts are involved, in order to illustrate the accuracy of the bounds obtained by Corollary \ref{Cor}, for $m=3$ and/or $m=4$.
\subsection{The Stieltjes-Wigert polynomials}
The Stieltjes-Wigert polynomials (cf.\ \cite[Section 14.27]{KLS})
 \[S_n(x;q)={1\over (q;q)_n}\qhypergeom{1}{1}{q^{-n}}{0}{q;-q^{n+1}x}\]
are orthogonal on $(0,\infty).$ When we let $m=3,k=0$ in \eqref{general4}, we obtain the equation
\begin{eqnarray}\label{BOUNDSW}
{q}^{6} S_{n-3} \left( x;q\right)&=&-q \left( {q}^{n}-1 \right)
   \left( x{q}^{2\,n}+{q}^{n+2}-{q}^{3}-{q}^{4} \right) S_{n} \left( x;q\right)\\
   &+&\nonumber \left( {q}^{4n}{x}^{2}+{q}^{2n+1}
   \left( q+1 \right)  \left( {q}^{n}-{q}^{2}-1 \right) x+
   {q}^{3} \left( {q}^{2\,n}+(q-{q}^{n})({q}^{2}+{q}+1)\right) \right)S_{n-1
   } \left( x;q\right).
 \end{eqnarray}
The coefficient of $S_{n-1}$ provides us with the two inner bounds
   \begin{equation} \nonumber
   B_3= \frac { ({q}^{2}-{q}^{n}+1)(q+1) \pm \sqrt {{q}^{6}-2\,{q}^{
   			n+4}+2\,{q}^{5}+  {q}^{2n+2} -{q}^{4}-2\,
   		{q}^{2n+1}+ {q}^{2n}-{q}^{2}-2\,
   		{q}^{n}+2\,q+1}  }{2  {q}^{2n-1} }.
   \end{equation}	
In Table \ref{SWbounds} we show the quality of the largest of these two bounds, which is a good lower bound for the largest zero of $S_n(x;q)$. A more accurate lower bound is obtained from
   \begin{eqnarray}\nonumber
   {q}^{12}S_{n-4}(\alpha,q)&=&q \left({q}^{n} -1\right)
   ( {q}^{4\,n}{x}^{2}+ ({q}^{3\,n+2}-{q}^{2\,n+5}-{q}^
{2\,n+3})(q+1) x\\
&+&\nonumber {q}^{2\,n+5}+({q}^{8}-{q}^{n+6})({q}^{2}+{q}+1)) S_{n}(\alpha,q) - A_3(x) S_{n-1}(\alpha,q),
\end{eqnarray}
with
\begin{eqnarray}A_3(x) \label{SW4}
&=&{q}^{6\,n}{x}^{3}+\left( {q}^{5\,n+1}-{q}^{
4\,n+4}-{q}^{4\,n+1}\right)\left({q}^{2}+{q}+1\right)  {x}^{2}\\
&+& \nonumber\left( \left({q}^{2\,n+8}+{q}^{2\,n+6}+{q}^{2\,n+4}+{q}^{4\,n+3}-{q}^{3\,n+6}-{q}^{3\,n+5}-{q}^{3\,n+4}\right)\left({q}^{2}+{q}+1\right)
-{q}^{3\,n+3}(q+1)\right ) x\\
&-&\nonumber{q}^{9}({q}^{3}+{q}^{2}+q+1)-{q}^{2\,n+6}({q}^{3}+{q}^{2}+q+1)+{q}^{3\,n+6}+{q}^{n+7}({q}^{2}+1)({q}^{2}+q+1).
\end{eqnarray}
The values of the zeros of $A_3(x)$ can be found numerically and we also show the value of the largest zero of $A_3(x)$ in Table \ref{SWbounds}.
\begin{table}[!ht]\centering  \caption{Examples to show the quality of the inner bound for the extreme zeros of  $S_n(x;q)$ for different values of $n$ and $q$.}\label{SWbounds}
   \begin{tabular}{|c|c|c|c|c|}
   \hline
   $n$&$q$&$B_3$ from (\ref{BOUNDSW})&\small{Lowest zero of $A_3$ (\ref{SW4})}& $x_{n,n}$\\
   \hline
   10&0.5&$8.3925988\cdot10^5$&$8.3946795\cdot10^5$&$8.3946799\cdot10^5$\\\hline
   10&0.9&15.3689887&16.0730951&16.1508699\\\hline
   70&0.5
   &$1.116350296\cdot10^{42}$&$1.1166280\cdot10^{42}$&$1.1166281\cdot10^{42}$\\\hline
   70&0.9
   &$5.9402911\cdot10^6$&$6.2658907\cdot10^6$&$6.3132591\cdot10^6$\\\hline
   \end{tabular}
   \end{table}		

  \subsection{The discrete $q$-Hermite II polynomials}
The discrete $q$-Hermite II polynomials (cf.\ \cite[Section 14.29]{KLS})
\[\tilde{h}_n(x;q)=i^{-n}q^{-\binom{n}{2}}\qhypergeom{2}{0}{q^{-n},ix}{-}{q;-q^n},\]
are orthogonal on the real line. The zeros of these polynomials are symmetric about the origin with a simple zero at the origin when $n$ is odd. We get the recurrence equation
\begin{equation}
\left( {q}^{-n+1};q \right)_4q^{14} \tilde{h}_{n-5}(x;q)=- \left({q}^{n+3}+{q}^{n+4}-{q}^{5}-{q}^{7}+{q}^{2\,n}{x}^{2} \right)  {q}^{2\,n}x \tilde{h}_{n}(x;q)
 + G_4(x) \tilde{h}_{n-1}(x;q)\nonumber
\end{equation}
with $$G_4(x)={q}^{4\,n}{x}^{4}+ \left({q}^{2}+{q}^{n}-q-{q}^{3} \right)  \left( {q}^{2}+q+1 \right) {x}^{2}{q}^{2\,n+2}+{q}^{6} \left( {q}^{n}-q
 \right)  \left( {q}^{n}-{q}^{3} \right).$$
 The largest zero of $G_4(x)$ is
\begin{eqnarray}&&B_5=\label{Bound}\\
&&\nonumber \Bigg(~\frac{( q+{q}^{3}-{q}^{2}-{q}^{n})({q}^{2}+q+1)
  +\sqrt{({q}^{2}+{q}^{n}-q-{q}^{3}) ^{2} ({q}^{2}+q+1)^{2}-4({q}^{n} -{q}^{3} )( {q}^{n}-q ) {q}^{2} }}{2q^{2n-2}}~\Bigg)^{\frac12}.
  \end{eqnarray}
We show the quality of this bound in Table \ref{qHbounds}.
\begin{table}[!ht]\centering \caption{Examples to show the quality of the inner bound for the extreme zeros of  $\tilde{h}_{n}(x;q)$ for different values of $n$ and $q$.}\label{qHbounds}
   \begin{tabular}{|c|c|c|c|}
   \hline
   $n$&$q$&$B_5$ in (\ref{Bound})& $x_{n,n}$\\
   \hline
   10&0.5&$406.3811$&$406.4079$\\\hline
   10&0.98&0.72968&0.76655\\\hline
   100&0.5&$5.0368\cdot10^{29}$&$5.0372\cdot10^{29}$\\\hline
   100&0.98 &$10.7735$&$12.1420$\\\hline
   \end{tabular}
   \end{table}	

\section{Acknowledgments}
This work has been supported by a TWAS / DFG fellowship for A. Jooste and  the Institute of Mathematics of the University of Kassel (Germany) for D.D. Tcheutia. All these institutions receive our sincere thanks.

\end{document}